\documentclass[final,letterpaper,11pt]{article}%
\usepackage[textwidth=5.5in,textheight=9in,centering,letterpaper]{geometry}
\usepackage{titletoc}
\usepackage{amsxtra,amscd}
\usepackage{graphicx}
\usepackage[amsmath,hyperref,thmmarks]{ntheorem}
\usepackage{natbib}
\usepackage{verbatim}
\usepackage{paralist}
\usepackage{makeidx}
\usepackage{amsmath}
\usepackage{amsfonts}
\usepackage{amssymb}
\usepackage{empheq}
\usepackage[scaled=0.92]{helvet}
\usepackage[notref,notcite]{showkeys}
\usepackage{url}%
\theoremnumbering{arabic}
\theoremheaderfont{\scshape}
\RequirePackage{latexsym}
\theorembodyfont{\slshape}
\theoremseparator{}
\newtheorem{X}{X}[section]

\newtheorem{conjecture}[X]{Conjecture}
\newtheorem{corollary}[X]{Corollary}

\newtheorem{lemma}[X]{Lemma}
\newtheorem{proposition}[X]{Proposition}
\newtheorem{theorem}[X]{Theorem}
\theorembodyfont{\upshape}

\newtheorem{definition}[X]{Definition}
\newtheorem{example}[X]{Example}

\newtheorem{remark}[X]{Remark}
\newtheorem{plain}[X]{}

\theorembodyfont{\small}
\newtheorem*{nt}{Notes}
\theorembodyfont{\normalsize}

\theoremstyle{nonumberplain}
\theoremheaderfont{\sc}
\theorembodyfont{\normalfont}
\theoremsymbol{\ensuremath{_\Box}}
\RequirePackage{amssymb}
\newtheorem{proof}{Proof.}
\qedsymbol{\ensuremath{_\Box}}
\theoremclass{LaTeX}
\makeindex
\newcommand\bquote{\begin{quote}}
\newcommand\equote{\end{quote}}
\newcommand\bsmall{\begin{small}}
\newcommand\esmall{\end{small}}

\let\cite=\citealt

\setdefaultitem{--\hspace{0.07in}}{}{}{}

\renewcommand{\Gamma}{\varGamma}
\renewcommand{\Pi}{\varPi}
\renewcommand{\Sigma}{\varSigma}

\def\1{{1\mkern-7mu1}}

\DeclareMathOperator{\Aut}{Aut}

\DeclareMathOperator{\End}{End}

\DeclareMathOperator{\Ext}{Ext}

\DeclareMathOperator{\Gal}{Gal}

\DeclareMathOperator{\Hom}{Hom}
\DeclareMathOperator{\id}{id}

\DeclareMathOperator{\im}{Im}  %Image \Im is already used by TeX.
\DeclareMathOperator{\Ind}{Ind}

\DeclareMathOperator{\inv}{inv}

\DeclareMathOperator{\Ker}{Ker}

\DeclareMathOperator{\ord}{ord}

\DeclareMathOperator{\rank}{rank}

\begin{document}

\title{Motivic complexes over finite fields and the ring of correspondences at the
generic point}
\author{James S. Milne
\and Niranjan Ramachandran}
\date{December 3, 2005: first version on web.\\
July 19, 2006: submitted version.\\
May 2, 2008: final version.}
\maketitle

\begin{abstract}
Already in the 1960s Grothendieck understood that one could obtain an almost
entirely satisfactory theory of motives over a finite field when one assumes
the full Tate conjecture. In this note we prove a similar result for motivic
complexes. In particular Beilinson's $\mathbb{Q}{}$-algebra of
\textquotedblleft correspondences at the generic point\textquotedblright\ is
then defined for all connected varieties. We compute it for all smooth
projective varieties (hence also for varieties birational to such a variety).

\end{abstract}
\tableofcontents

\subsubsection{Introduction}

More than forty years after Grothendieck predicted that the standard
cohomology functors factor through a tannakian category of pure motives, we
still do not know how to construct such a category. However, when the field is
finite and one assumes the full Tate conjecture, there is an almost entirely
satisfactory theory of pure motives. According to Deligne
(1994\nocite{deligne1994}, 1.4), this was known to Grothendieck, but it was
re-discovered by \citet{langlandsR1987}, who used it to state a conjecture,
more precise than earlier attempts by Langlands, on the structure of the
points modulo a prime on a Shimura variety. For a detailed description of the
category, see \cite{milne1994}.

It is generally hoped that the standard cohomology functors to triangulated
categories will factor through a triangulated category of motivic complexes
with $t$-structure whose heart is (defined to be) the category of mixed
motives (see, for example, \cite{deligne1994}, \S 3). We show that, over a
finite field, a triangulated category of motivic complexes exists with the
expected properties if and only if the Tate conjecture holds and homological
equivalence coincides with rational equivalence with $\mathbb{Q}{}%
$-coefficients (see Theorems \ref{03} and \ref{06} for more precise
statements). Moreover, then a category of effective motivic complexes exists
with the properties (A,B,C) of \cite{beilinson2002}, and so there is a
well-defined semisimple $\mathbb{Q}{}$-algebra of \textquotedblleft
correspondences at the generic point\textquotedblright\ attached to every
variety over a finite field. We compute this $\mathbb{Q}{}$-algebra for smooth
projective varieties (hence also for varieties birational to such a variety).
As this requires the generalized Tate conjecture (in the sense of
\cite{grothendieck1968}, \S 10), we begin by giving an elementary proof that
this follows from the usual Tate conjecture.

\subsubsection{Notations\label{notations}}

A variety is a geometrically-reduced separated scheme of finite type over a
field. For a variety $X$ over a perfect field $k$ of characteristic $p\neq0$
and algebraic closure $\bar{k}$, we set%
\begin{align*}
H_{l}^{i}(X)  &  =H_{\mathrm{et}}^{i}(X_{\bar{k}},\mathbb{Q}{}_{l}%
),\quad\text{if }l\neq p,\text{ and}\\
H_{p}^{i}(X)  &  =H_{\mathrm{crys}}^{i}(X/W)\otimes\mathbb{Q}{},\quad W=W(k).
\end{align*}
We use $(r)$ to denote a Tate twist, and we write $\mathrm{hom}(l)$ for the
equivalence relation on the space $Z^{\ast}(X)$ of algebraic cycles defined by
$H_{l}$. Similarly, we write $\mathrm{num}$ and $\mathrm{rat}$ for numerical
and rational equivalence. For an adequate equivalence relation $\sim$,
$Z_{\sim}^{i}(X)=Z^{i}(X)/\!\sim$ and $Z_{\sim}^{i}(X)_{\mathbb{Q}{}}=Z_{\sim
}^{i}(X)\otimes\mathbb{Q}{}$. For example, $Z_{\text{rat}}^{i}(X)$ is the Chow
group $CH^{i}(X)$.

By a functor between additive categories, we mean an additive functor. A
functor $F\colon\mathcal{C}{}\rightarrow\mathcal{C}^{\prime}$ of triangulated
categories together with an isomorphism of functors $F\circ T\approx
T^{\prime}\circ F$ is said to be \emph{triangulated} (formerly, exact;
\cite{verdier1977}, p4) if it takes distinguished triangles to distinguished triangles.

A triangulated category with $t$-structure (\cite{gelfandM1996}, IV 4.2, p278)
will be referred to simply as a $t$-category. All $t$-structures will be
assumed to be bounded\emph{ }(i.e., $\bigcup_{n\geq0}\mathcal{D}{}^{\leq
n}=\mathcal{D}{}=\bigcup_{n\geq0}\mathcal{D}{}^{\geq-n}$) and
nondegenerate\emph{ }(i.e., $\bigcap_{n\geq0}\mathcal{D}{}^{\leq-n}%
=0=\bigcap_{n\geq0}\mathcal{D}{}^{\geq n}$).

The symbol $\mathbb{F}{}$ denotes an algebraic closure of $\mathbb{F}{}_{p}$,
and the algebraic closure of $\mathbb{Q}{}$ in $\mathbb{C}{}$ is denoted
$\mathbb{Q}{}^{\mathrm{al}}$. Reductive groups are not required to be
connected. Isomorphisms are denoted $\approx$ and canonical isomorphisms
$\simeq$.

\section{The generalized Tate conjecture}

In this section, $k$ is the subfield $\mathbb{F}{}_{q}$ of $\mathbb{F}{}$, and
$l\neq p$.

\begin{plain}
\label{gtc0}By the \emph{full Tate conjecture} for a smooth complete variety
$X$ over $k$ and an $r\geq0$, we mean the statement that the order of the pole
of the zeta function $Z(X,t)$ at $t=q^{-r}$ is equal to the rank of the group
of numerical equivalence classes of algebraic cycles of codimension $r$ on
$X$. If the full Tate conjecture holds for $X$ and $r$, then, for all $l\neq
p$,

\begin{enumerate}
\item[$T^{r}(X,l)$:] the cycle class map ${}Z^{r}(X)\otimes\mathbb{Q}{}%
_{l}\rightarrow H_{l}^{2r}(X)(r)^{\Gal(\mathbb{F}{}/k)}$ is surjective, and

\item[$E^{r}(X,l)$:] the quotient map $Z_{\mathrm{hom}(l)}^{r}(X)_{\mathbb{Q}%
{}}\rightarrow Z_{\mathrm{num}}^{r}(X)_{\mathbb{Q}{}}$ is injective (i.e.,
$\mathrm{hom}(l)$ and $\mathrm{num}$ coincide with $\mathbb{Q}{}$-coefficients).
\end{enumerate}

\noindent Conversely, if $T^{r}(X,l)$ and $E^{r}(X,l)$ hold for a single $l$,
then the full Tate conjecture holds for $X$ and $r$ (\cite{tate1994}, \S 2).
The statement $T^{r}(X,l)$ is the \emph{Tate conjecture for }$X$, $r$, and
$l$.\footnote{More precisely, it is Conjecture 1 of \cite{tate1965}. Statement
$E^{r}(X,l)$ is a variant of the \textquotedblleft conjectural
statement\textquotedblright\ (\textbf{a}$^{\prime}$) of \cite{tate1965}. Our
notation follows that of \cite{tate1994}.}
\end{plain}

\subsection{Statement of the generalized Tate conjecture}

Define a \emph{Tate structure }to be a finite-dimensional $\mathbb{Q}{}_{l}%
$-vector space with a linear (Frobenius) map $\varpi$ whose characteristic
polynomial lies in $\mathbb{Q}{}[T]$ and whose eigenvalues are Weil
$q$-numbers, i.e., algebraic numbers $\alpha$ such that, for some integer $m$
(called the weight of $\alpha$), $\left\vert \rho(\alpha)\right\vert =q^{m/2}$
for every homomorphism $\rho\colon\mathbb{Q}{}[\alpha]\rightarrow\mathbb{C}$,
and, for some integer $n$, $q^{n}\alpha$ is an algebraic integer.
\noindent\ When the eigenvalues are all of weight $m$ (resp. algebraic
integers, resp. semisimple), we say that $V$ is of \emph{weight} $m$ (resp.
\emph{effective}, resp. \emph{semisimple}). For example, for any smooth
complete variety $X$ over $k$, $H_{l}^{i}(X)$ is an effective Tate structure
of weight $i/2$ (\cite{deligne1980}), which is semisimple if $X$ is an abelian
variety (\cite{weil1948}, no. 70) or if the full Tate conjecture holds for
$X\times X$ (\cite{milne1986v}, 8.6).

Let $X$ be a smooth complete variety over $k$. For each $r$, let $F_{a}%
^{r}H_{l}^{i}(X)\subset H_{l}^{i}(X)$ denote the subspace of classes with
support in codimension at least $r$, i.e.,%
\[
F_{a}^{r}H_{l}^{i}(X)=\bigcup\nolimits_{U}\Ker(H_{l}^{i}(X)\rightarrow
H_{l}^{i}(U))
\]
where $U$ runs over the open subvarieties of $X$ such that $X\smallsetminus U$
has codimension $\geq r$.

\begin{example}
\label{gtc1a}If $Z$ is a \textit{smooth }closed subvariety of $X$ of
codimension $r$, then there is an exact Gysin sequence%
\[
\cdots\rightarrow H_{l}^{i-2r}(Z)(-r)\rightarrow H_{l}^{i}(X)\rightarrow
H_{l}^{i}(U)\rightarrow\cdots,\quad U=X\smallsetminus Z,
\]
(e.g., \cite{milne1980}, VI 5.4), and so the kernel of $H_{l}^{i}%
(X)\rightarrow H_{l}^{i}(U)$ is an effective Tate structure of weight $i$
whose twist by $\mathbb{Q}{}_{l}(r)$ is still effective.
\end{example}

\begin{conjecture}
\label{gtc1}(Generalized Tate conjecture; cf. \cite{grothendieck1968}, 10.3.).
For a smooth complete variety $X$ over $k$, every semisimple Tate substructure
$V\subset H_{l}^{i}(X)$ such that $V(r)$ is still effective is contained in
$F_{a}^{r}H_{l}^{i}(X)$.
\end{conjecture}

\begin{remark}
\label{gtc2}Let $X$ be a smooth complete variety over $k$. For any $i$ and
$r$, the set of eigenvalues $\alpha$ of $\varpi_{X}$ on $H_{l}^{i}(X)$ such
that $\alpha/q^{r}$ is an algebraic integer is stable under Galois
conjugation. Therefore, there is a subspace $F_{b}^{r}H_{l}^{i}(X)$ of
$H_{l}^{i}(X)$ that becomes the sum of the eigenspaces of these $\alpha$ over
$\mathbb{Q}{}_{l}^{\mathrm{al}}$. It is the largest semisimple Tate
substructure of $H_{l}^{i}(X)$ whose twist by $\mathbb{Q}{}_{l}(r)$ is still
effective, and so the generalized Tate conjecture \ref{gtc1} is the statement:%
\[
F_{b}^{r}H_{l}^{i}(X)\subset F_{a}^{r}H_{l}^{i}(X).
\]

\end{remark}

\begin{example}
\label{gtc3}Let $Z^{\prime}$ be a closed irreducible subvariety of
$X_{\mathbb{F}{}}$ of codimension $r$. Then%
\begin{equation}
H_{Z^{\prime}}^{2r}(X_{\mathbb{F}{}},\mathbb{Q}{}_{l}(r))\rightarrow
H^{2r}(X_{\mathbb{F}{}},\mathbb{Q}{}_{l}(r))\rightarrow H^{2r}(X_{\mathbb{F}%
{}}\smallsetminus Z^{\prime},\mathbb{Q}{}_{l}(r)) \label{eq1}%
\end{equation}
is exact, and $H_{Z^{\prime}}^{2r}(X_{\mathbb{F}{}},\mathbb{Q}{}_{l}%
(r))\simeq\mathbb{Q}{}_{l}$; moreover, the image of $1$ under the first map is
the cohomology class of $Z^{\prime}$ (cf. \cite{milne1980}, p269). For any
open $U\subset X$, the kernel of
\[
H_{l}^{2r}(X)(r)\rightarrow H_{l}^{2r}(U)(r)
\]
is spanned by the cohomology classes of the irreducible components of
$(X\smallsetminus U)_{\mathbb{F}{}}$, and some power of $\varpi_{X}$ acts as
$1$ on it. On the other hand, $F_{b}^{r}H_{l}^{2r}(X)(r)$ is the largest
subspace of $H_{l}^{2r}(X)(r)$ on which some power of $\varpi$ acts as $1$.
Thus, the generalized Tate conjecture with $i=2r$ states that this subspace is
spanned by the classes of algebraic cycles of codimension $r$ on
$X_{\mathbb{F}{}}$. This is the Tate conjecture stated over $\mathbb{F}{}$
rather than $\mathbb{F}{}_{q}$.
\end{example}

\subsection{The Tate conjecture implies the generalized Tate conjecture}

Recall that, for a proper map $\pi\colon Y\rightarrow X$ of smooth varieties
over an algebraically closed field, the Gysin map%
\[
\pi_{\ast}\colon H^{i}(Y,\mathbb{Q}{}_{l})\rightarrow H^{i-2c}(X,\mathbb{Q}%
{}_{l}(-c)),\quad c=\dim Y-\dim X,
\]
is defined to be the Poincar\'{e} dual of
\[
\pi^{\ast}\colon H_{c}^{2d-i}(X,\mathbb{Q}{}_{l}(d))\rightarrow H_{c}%
^{2d-i}(Y,\mathbb{Q}{}_{l}(d)),\quad d=\dim Y
\]
(\cite{milne1980}, VI 11.6). We shall need to know that these maps are
compatible with restriction to open subvarieties.

\begin{lemma}
\label{gtc7}Let $\pi\colon Y\rightarrow X$ be a proper map of smooth complete
varieties over an algebraically closed field, and let $j\colon
U\hookrightarrow X$ an open immersion. Then the commutative diagram at left
gives rise to the commutative diagram at right:%
\[
\begin{CD}
Y @<{j^{\prime}}<< \pi^{-1}U\\
@VV{\pi}V@VV{\pi^{\prime}}V\\
X @<j<< U
\end{CD}\qquad\begin{CD}
H^{i}(Y,\mathbb{Q}_{l}) @>{j^{\prime\ast}}>> H^{i}(\pi^{-1}U,\mathbb{Q}_{l})\\
@VV{\pi_{\ast}}V@VV{\pi_{\ast}^{\prime}}V\\
H^{i-2c}(X,\mathbb{Q}_{l}(-c))@>{j^{\ast}}>>H^{i-2c}(U,\mathbb{Q}_{l}(-c))
\end{CD}
\]

\end{lemma}

\begin{proof}
Exercise for the reader.
\end{proof}

\begin{proposition}
\label{gtc5}Every effective semisimple Tate structure is isomorphic to a Tate
substructure of $H_{l}^{\ast}(A)$ for some abelian variety $A$ over
$\mathbb{F}{}_{q}$.
\end{proposition}

\begin{proof}
We may assume that the Tate structure $V$ is simple. Then $V$ has weight $m$
for some $m\geq0$, and the characteristic polynomial $P(T)$ of $\varpi$ is a
monic irreducible polynomial with coefficients in $\mathbb{Z}{}$ whose roots
all have real absolute value $q^{m/2}$. According to Honda's theorem
(\cite{honda1968}; \cite{tate1968}), $P(T)$ is the characteristic polynomial
of an abelian variety $A$ over $\mathbb{F}{}_{q^{m}}$. Let $B$ be the abelian
variety over $\mathbb{F}{}_{q}$ obtained from $A$ by restriction of the base
field. The eigenvalues of the Frobenius map on $H_{l}^{1}(B)$ are the
$m^{\mathrm{th}}$-roots of the eigenvalues of the Frobenius map on $H_{l}%
^{1}(A)$, and it follows that $V$ is a Tate substructure of $H_{l}^{m}(B)$.
\end{proof}

\begin{lemma}
\label{gtc6}Let $z$ be an algebraic cycle of codimension $\dim T{}{}+r$ on the
product $T{}{}\times X$ of two smooth complete varieties over $k$ (i.e., $z$
is an algebraic correspondence of degree $r$ from $T{}{}$ to $X$). Assume that
the push-forward of $z$ on $X$ is nonzero. Then the image of the map
\[
z_{\ast}\colon H_{l}^{i-2r}(T{}{})(-r)\rightarrow H_{l}^{i}(X)
\]
defined by $z$ is contained in $F_{a}^{r}H_{l}^{i}(X)$.
\end{lemma}

\begin{proof}
Let $p,q$ denote the projection maps $T{}{}\times X\rightrightarrows T{}{},X$,
and let $[z]$ denote the cohomology class of $z$ in $H_{l}^{2d_{T{}{}}%
+2r}(T\times X)(d_{T{}{}}+r)$, $d_{T{}{}}=\dim T{}{}$. Then
\[
z_{\ast}(a)\overset{\text{{\tiny def}}}{=}q_{\ast}([z]\cup p^{\ast}(a)),\quad
a\in H_{l}^{i-2r}(T)(-r).
\]
As the push-forward $q_{\ast}(z)$ of $z$ is nonzero, its support $Z$ has
codimension $r$.\footnote{Recall that the push-forward $q_{\ast}(z)$ of an
irreducible $z$ is defined to be zero if $\dim(q(z))<\dim z$.} Let
$U=X\smallsetminus Z$. Then $z$ has support in $T{}{}\times Z$, and so $[z]$
maps to zero in $H_{l}^{2d_{T{}{}}+2r}(T\times U)(d_{T}+r)$. According to
(\ref{gtc7}), the diagram%
\[
\begin{CD}
H_l^{i+2d_{T}}(T\times X)(d_T)
@>>> H_l^{i+2d_{T}}(T\times U)(d_{T})\\
@VV{q_{\ast}}V@VV{q_{\ast}}V\\
H_l^{i}(X) @>>> H_l^{i}(U)
\end{CD}
\]
commutes, which shows that $z_{\ast}(a)$ maps to zero in $H_{l}^{i}(U)$, and
therefore lies in $F_{a}^{r}H_{l}^{i}(X)$.
\end{proof}

\begin{lemma}
\label{gtc6l}Let $X$ be a smooth complete variety over $k$ and let
$i,r\in\mathbb{N}{}$. If there exists a smooth complete variety $T$ such that

\begin{itemize}
\item $H_{l}^{i-2r}(T)$ is a semisimple Tate structure,

\item the Tate conjecture $T^{\dim(T)+r}(T\times X,l)$ holds, and

\item $F_{b}^{r}H_{l}^{i}(X)(r)$ is isomorphic to a Tate substructure of
$H_{l}^{i-2r}(T)$
\end{itemize}

\noindent then $F_{b}^{r}H_{l}^{i}(X)\subset F_{a}^{r}H_{l}^{i}(X)$.
\end{lemma}

\begin{proof}
Let $d=\dim(T)$ and let $V$ be a Tate substructure of $H_{l}^{i-2r}(T)$ for
which there exists an isomorphism $f\colon V(-r)\rightarrow F_{b}^{r}H_{l}%
^{i}(X)$. Then%
\begin{align*}
H_{l}^{2d+2r}(T\times X)(d+r)  &  \supset H_{l}^{2d+2r-i}(T)(d+r)\otimes
H_{l}^{i}(X)\\
&  \simeq\Hom(H_{l}^{i-2r}(T)(-r),H_{l}^{i}(X))\\
&  \supset\Hom(V(-r),F_{b}^{r}H_{l}^{i}(X))\ni f.
\end{align*}
(The last inclusion depends on the choice of stable complement for $V$ in
$H_{l}^{i-2r}(T)$.) As $f$ is fixed by $\Gal(\mathbb{F}{}/k)$, it can be
approximated by the cohomology class of an algebraic correspondence $z$ of
degree $r$ from $T$ to $X$. Moreover, $z$ can be chosen so that $z_{\ast}$ is
injective on $V$. Obviously $z_{\ast}$ maps $H_{l}^{i-2r}(T)(-r)$ into
$F_{b}^{r}H_{l}^{i}(X)$, and so
\[
F_{b}^{r}H_{l}^{i}(X)\subset z_{\ast}V(-r)\overset{\ref{gtc6}}{\subset}%
F_{a}^{r}H_{l}^{i}(X).
\]

\end{proof}

\begin{theorem}
\label{gtc8}Let $X$ be a smooth complete variety over $k$. If the Tate
conjecture holds for all varieties of the form $A\times X$ with $A$ an abelian
variety (and some $l$), then the generalized Tate conjecture holds for $X$
(and the same $l$).
\end{theorem}

\begin{proof}
As we noted above, $H_{l}^{\ast}(A)$ is a semisimple, and so this follows from
(\ref{gtc5}) and \ref{gtc6l}).
\end{proof}

\begin{corollary}
\label{gtc9}If the Tate conjecture holds for all abelian varieties over $k$
(or for all smooth complete varieties over $k$) and some $l$, then the
generalized Tate conjecture holds for the same class and that $l$.
\end{corollary}

\begin{remark}
\label{gtc8a}As others have noted (\cite{kahn2002}, Theorem 2;
\cite{andre2004}, 8.2), when one assumes the full Tate conjecture, the
generalized Tate conjecture follows directly from the description of the
simple motives in terms of Weil numbers (see \cite{milne1994}, Proposition 2.6).
\end{remark}

\subsection{Complements}

\begin{plain}
\label{gtc9a}Let $X$ be a smooth projective variety over $k$, and let
$V=F_{b}^{r}H_{l}^{i}(X)$. We know that $V(-r)\subset H_{l}^{i-2r}(A)$ for
some abelian variety $A$ over $k$ (see \ref{gtc5}). If $\dim A=d>i-2r$, then,
according to the Lefschetz hypersurface-section theorem, for any smooth
hypersurface section $Y$ of $A$ (which exists by \cite{gabber2001}),
$V(-r)\subset H_{l}^{i-2r}(Y)$. Continuing in this fashion, we get that
$V(-r)\subset H_{l}^{i-2r}(T)$ for some smooth projective $T$ of dimension
$i-2r$. Therefore, under the assumption of the Tate conjecture, there exists a
smooth projective variety $T$ of dimension at most $i-2r$ over $k$ and an
algebraic correspondence $z$ from $T$ to $X$ of degree $r$ such that $z_{\ast
}H_{l}^{i-2r}(T)(r)=$ $F_{b}^{r}H_{l}^{i}(X)$.
\end{plain}

\begin{plain}
\label{gtc2a} \citet[8.2.8]{deligne1974t} proves the following:\bquote Let $X$
be a smooth complete variety over $\mathbb{C}{}$, and let $Z$ be a closed
subvariety of $X$ of codimension $r$. For any desingularization $\tilde
{Z}\rightarrow Z$ of $Z$, the sequence%
\[
H^{i-2r}(\tilde{Z},\mathbb{Q}{}(-r))\rightarrow H^{i}(X,\mathbb{Q}%
{})\rightarrow H^{i}(U,\mathbb{Q}),\quad U={}X\smallsetminus Z,
\]
is exact.\equote A similar argument\footnote{For any proper surjective
morphism $f\colon Y\rightarrow Z$ from a smooth projective variety $Y$, we can
find a smooth projective simplicial scheme $Y_{\bullet}$ with $Y_{0}=Y$ that
is a proper hypercovering of $Z$. The corresponding spectral sequence
($l$-adic analogue of the spectral sequence \cite{deligne1974t}, 8.1.19.1)
degenerates at $E_{2}$ with $\mathbb{Q}{}_{l}$-coefficients because of weight
considerations, and gives an exact sequence%
\[
0\rightarrow\frac{H_{l}^{i}(Z)}{W_{i-1}H_{l}^{i}(Z)}\rightarrow H_{l}%
^{i}(Y_{0})\overset{\delta_{0}-\delta_{1}}{\longrightarrow}H_{l}^{i}(Y_{1}).
\]
It follows that the image of $H_{l}^{i}(Z)$ in $H_{l}^{i}(Y)$ is the (largest)
quotient of pure weight $i$ of $H_{l}^{i}(Z)$. This implies the $l$-adic
analogue of \cite{deligne1974t}, 8.2.7, (the proof there works as the
$Gr_{\ast}^{W}$ functor is exact) and of ibid. 8.2.8.} proves the following
$l$-adic analogue: \bquote Let $X$ be a smooth complete variety over a perfect
field $k{}$, and let $Z$ be a closed subvariety of $X$ of codimension $r$. For
any smooth alteration $\tilde{Z}\rightarrow Z$ of $Z$, the sequence%
\[
H_{l}^{i-2r}(\tilde{Z})(-r)\rightarrow H_{l}^{i}(X)\rightarrow H_{l}%
^{i}(U),\quad U=X\smallsetminus Z,
\]
is exact.\equote Since \citet[3.1]{dejong1996} shows that smooth alterations
always exist, this implies that%
\[
F_{a}^{r}H_{l}^{i}(X)\subset F_{b}^{r}H_{l}^{i}(X).
\]
The generalized Tate conjecture then states that
\[
F_{a}^{r}H_{l}^{i}(X)=F_{b}^{r}H_{l}^{i}(X).
\]

\end{plain}

\begin{plain}
\label{gtc9b}The above statements hold \textit{mutatis mutandis }for $p$. For
a smooth complete variety $X$, $H_{p}^{i}(X)$ is an $F$-isocrystal, i.e., a
finite-dimensional vector space over $B(\mathbb{F}{}_{q})\overset
{\text{{\tiny def}}}{=}W(\mathbb{F}{}_{q})\otimes\mathbb{Q}{}$ equipped with a
$\sigma$-linear bijection $F\colon H_{p}^{i}(X)\rightarrow H_{p}^{i}(X)$. The
full Tate conjecture for $X$ and $r$ is equivalent to

\begin{enumerate}
\item[$T^{r}(X,p)$:] the cycle class map $Z^{r}(X)\otimes\mathbb{Q}{}%
_{p}\rightarrow H_{p}^{2r}(X)(r)^{F=1}$ is surjective (\emph{Tate conjecture
for }$p$), and

\item[$E^{r}(X,p)$:] the quotient map $Z_{\mathrm{hom}(p)}^{r}(X)_{\mathbb{Q}%
{}}\rightarrow Z_{\mathrm{num}}^{r}(X)_{\mathbb{Q}{}}$ is injective
\end{enumerate}

\noindent(cf. \cite{milne2007rtc}, \S 1). \noindent Define%
\[
F_{a}^{r}H_{p}^{i}(X)=\bigcup\nolimits_{Z}\im(H_{p}^{i-2r}(\tilde
{Z})(-r)\rightarrow H_{p}^{i}(X))
\]
where $Z$ runs over the closed subvarieties of $X$ such that $Z$ is of
codimension at least $r$ and $\tilde{Z}$ is a smooth alteration of $Z$. If the
Tate conjecture holds for smooth complete varieties over $k$ and $p$, then%
\[
F_{a}^{r}H_{p}^{i}(X)=F_{b}^{r}H_{p}^{i}(X)
\]
where $F_{b}^{r}H_{p}^{i}(X)\subset H_{p}^{i}(X)_{[r,\infty)}$ is the largest
semisimple sub-isocrystal of $H_{p}^{i}(X)$ with slopes at least $r$. The
proofs are similar to those in the case $l\neq p$ --- we omit the details.
\end{plain}

\begin{plain}
\label{gtc9c} Similar arguments show that the generalized Tate conjecture over
number fields follows from the Tate conjecture and an effective version of the
Fontaine-Mazur conjecture (\cite{fontaineM1995}, Conjecture 1, p44) that
specifies which representations arise from effective motives.
\end{plain}

\begin{nt}
It was known to Grothendieck that the generalized Hodge conjecture follows
from the usual Hodge conjecture and the following weak analogue of
(\ref{gtc5}), \bquote Let $V$ be a simple Hodge substructure of the cohomology
of a smooth complex projective variety; if its Tate twist $V(r)$ is still
effective (i.e., has only nonnegative Hodge numbers), then $V(r)$ occurs in
the cohomology of a smooth complex projective variety.\equote presumably by
more-or-less the above argument. See \cite{grothendieck1969h}, top of p301
(also \cite{schoen1989}, \S 0).
\end{nt}

\section{The category of pure motives}

In this section $k=\mathbb{F}{}_{q}$.

For any adequate equivalence relation $\sim$, Grothendieck's construction
gives a rigid pseudo-abelian tensor $\mathbb{Q}{}$-category $\mathcal{M}%
{}_{\sim}(k)$ of pure motives (\cite{saavedra1972}, VI 4.1.3.5, p359) and a
map $h$ from the smooth projective varieties over $k$ to $\mathcal{M}{}_{\sim
}(k)$ which is natural for algebraic correspondences modulo $\sim$. Because
rational equivalence is the finest adequate equivalence relation, $h$ factors
through a tensor functor $\mathcal{M}{}_{\text{rat}}(k)\rightarrow
\mathcal{M}{}_{\sim}(k)$. Conversely, a tensor functor from $\mathcal{M}%
{}_{\text{rat}}(k)$ to an additive tensor category with $\End(\1)=\mathbb{Q}%
{}$ defines an adequate equivalence relation (cf. \cite{jannsen2000}, 1.7).
When $\sim$ is numerical equivalence, $\mathcal{M}_{\sim}{}(k{})$ is a
semisimple (\cite{jannsen1992}).

For a smooth projective variety $X$ over $k$, there are well-defined
polynomials $P_{X,i}(T)\in\mathbb{Q}{}[T]$ such that $P_{X,i}(T)=\det
(1-\varpi_{X}T\mid H_{l}^{i}(X))$ for all $l$; moreover, $P_{X,i}$ has
reciprocal roots of absolute value $q^{\frac{i}{2}}$ (\cite{deligne1974}). The
$P_{X,i}(T)$ are relatively prime, and so there exist $P^{i}(T)\in\mathbb{Q}%
{}[T]$, well-defined up to a multiple of $\prod_{i}P_{X,i}(T)$, such that%
\begin{equation}
P^{i}(T)\equiv\left\{
\begin{array}
[c]{lll}%
1 & \quad & \text{mod }P_{X,i}(T)\\
0 &  & \text{mod }P_{X,j}(T)\text{ for }j\neq i.
\end{array}
\right.  \label{e1}%
\end{equation}
Because $\prod_{i}P_{X,i}(\varpi_{X})$ acts as zero on $H_{l}^{\ast}(X)$, the
graph $p^{i}$ of $P^{i}(\varpi_{X})$ is a well-defined element of
$Z_{\text{hom}(l)}(X\times X)_{\mathbb{Q}{}}$ (or $Z_{\text{num}}(X\times
X)_{\mathbb{Q}{}}$), and $\{p^{0},\ldots,p^{2d}\}$ is a complete set of
orthogonal idempotents. Let $hX=\bigoplus\nolimits_{i}h^{i}X$ be the
corresponding decomposition. When we use this decomposition to modify the
commutativity constraint in $\mathcal{M}_{\text{num}}(k)$, the rank of each
object of $\mathcal{M}_{\text{num}}(k)$ becomes a nonnegative integer, and so
$\mathcal{M}_{\text{num}}{}(k{})$ is a tannakian category (\cite{deligne1990}, 7.1).

The category $\mathcal{M}_{\text{num}}(k)$ has a canonical (Frobenius) element
$\varpi\in\Aut^{\otimes}(\id_{\mathcal{M}_{\text{num}}(k)})$ and a canonical
(weight) $\mathbb{Z}{}$-gradation. An object $M$ of $\mathcal{M}_{\text{num}%
}(k)$ is of pure weight $m$ if and only if its Frobenius element $\varpi_{M}$
has eigenvalues of absolute value $q^{m/2}$.

Recall (\cite{deligne1989}, \S 6) that the fundamental group $\pi
(\mathcal{T}{})$ of a tannakian category is an affine group scheme in
$\Ind\mathcal{T}{}$ that acts on each object of $\mathcal{T}{}$ in such a way
that these actions define an isomorphism%
\[
\omega(\pi(\mathcal{T}{}))\simeq\underline{\Aut}^{\otimes}(\omega)
\]
for each fibre functor $\omega$. Any subgroup of the centre of $\pi
(\mathcal{T}{})$ lies in $\Ind\mathcal{T}{}^{0}$ where $\mathcal{T}{}^{0}$ is
the full subcategory of trivial objects (those isomorphic to a multiple of
$\1$). Since $\Hom_{\mathcal{T}{}}(\1,-)$ defines an equivalence of
$\mathcal{T}{}^{0}$ with the finite-dimensional vector spaces over the ground
field, such a subgroup can be identified with an affine group scheme in the
usual sense. For example, the centre of $\pi(\mathcal{T}{})$ is $\underline
{\Aut}^{\otimes}(\id_{\mathcal{T}{}}$) (cf. \cite{saavedra1972}, II 3.3.3.2).

Recall (e.g., \cite{milne1994}, \S 2) that the Weil-number group $P$ is the
affine group scheme of multiplicative type over $\mathbb{Q}{}$ whose character
group consists of the Weil $q$-numbers in $\mathbb{Q}{}^{\mathrm{al}}$. Define
the Frobenius element $\varpi_{\text{univ}}\ $in $P(\mathbb{Q}{})$ to be that
corresponding to $\alpha\mapsto\alpha$ under the bijection
\[
P(\mathbb{Q}{})\simeq\Hom(X^{\ast}(P),\mathbb{Q}{}^{\mathrm{al}}%
)^{\Gal(\mathbb{Q}{}^{\mathrm{al}}/\mathbb{Q}{})}.
\]
Note that, for any smooth projective variety $X$ over $\mathbb{F}{}_{q}$, the
roots of $P_{X,i}(T)$ in $\mathbb{Q}{}^{\mathrm{al}}$ are Weil $q$-integers of
weight $i$ (i.e., Weil $q$-numbers of weight $i$ that are algebraic integers).

\begin{lemma}
\label{0}The group of Weil $q$-numbers is generated by the Weil $q$-numbers of
abelian varieties over $k$.
\end{lemma}

\begin{proof}
Let $\alpha$ be a Weil $q$-number. After multiplying $\alpha$ by a power of
$q$, we may suppose that it is a Weil $q$-integer, of weight $m$ say. Then
$\alpha^{1/m}$ is a Weil $q$-integer of weight $1$, and hence arises from an
abelian variety by \cite{honda1968}.
\end{proof}

\begin{proposition}
\label{00}The affine subgroup scheme of $\pi(\mathcal{M}_{\text{num}}(k))$
generated by $\varpi_{\text{univ}}$ is canonically isomorphic to $P$. It
equals $\pi(\mathcal{M}_{\text{num}}(k))$ if and only if the full Tate
conjecture holds over $k$.
\end{proposition}

\begin{proof}
Let $Z=\underline{\Aut}^{\otimes}(\id)$ be the centre of $\pi(\mathcal{M}%
_{\text{num}}{}(k))$. Because $\mathcal{M}_{\text{num}}(k)$ is semisimple,
$\pi(\mathcal{M}_{\text{num}}(k))$ is pro-reductive (cf. \cite{deligneM1982},
2.23). Therefore $Z$ is of multiplicative type, which implies that the closed
subgroup scheme $\langle\varpi_{\text{univ}}\rangle$ generated by
$\varpi_{\text{univ}}$ is also of multiplicative type. The homomorphism
$\chi\mapsto\chi(\varpi_{\text{univ}})\colon X^{\ast}(\langle\varpi
_{\text{univ}}\rangle)\rightarrow\mathbb{Q}^{\mathrm{al}\times}$ is injective,
and its image consists of the Weil $q$-numbers that occur as roots of the
characteristic polynomial of $\varpi_{M}$ for some $M$ in $\mathcal{M}%
_{\text{num}}(k)$. According to Lemma \ref{0}, this consists of all Weil
$q$-numbers, and so $X^{\ast}(\langle\varpi_{\text{univ}}\rangle)\simeq
X^{\ast}(P)$. Hence $\langle\varpi_{\text{univ}}\rangle\simeq P$.

If the full Tate conjecture holds, then, for any fibre functor $\omega$ over
$\mathbb{Q}{}^{\mathrm{al}}$ and smooth projective variety $X$, the
$\mathbb{Q}{}^{\mathrm{al}}$-span of the algebraic cycles in $\omega
(h^{2i}(X)(i))$ consists of the tensors fixed by $\varpi_{\text{univ}}$.
Therefore, the inclusion $\langle\varpi_{\text{univ}}\rangle\hookrightarrow
\underline{\Aut}^{\otimes}(\omega)$ is an isomorphism, i.e., $\omega
(P)\hookrightarrow\omega(\pi(\mathcal{M}_{\text{num}}(k))$ is an isomorphism,
which implies that $P\hookrightarrow$ $\pi(\mathcal{M}_{\text{num}}(k))$ is an
isomorphism. The converse can be proved by the same argument as in the proof
of \cite{milne1999lm}, Proposition 7.4.
\end{proof}

If $\mathrm{num}$ and $\mathrm{hom}(l)$ coincide with $\mathbb{Q}{}%
$-coefficients, then $H_{l}$ defines a fibre functor $\omega_{l}$ on
$\mathcal{M}_{\text{num}}(k)$. Without any assumptions, it is known that there
exists a polarizable semisimple tannakian category with fundamental group $P$
and with fibre functors $\omega_{l}$ for all $l$. Moreover, any two such
systems are equivalent (\cite{langlandsR1987}; \cite{milne2003}, \S 6).
However, it has not been shown that there exists a natural functor from
$\mathcal{M}{}_{\text{rat}}(k)$ to such category. In fact, we have the following:

\begin{proposition}
\label{01}If there exists a full tensor functor $r$ preserving Frobenius
elements from $\mathcal{M}_{\text{rat}}(k)$ to a tannakian category
$\mathcal{M}{}$ with fundamental group $P$, then the full Tate conjecture
holds over $k$, and $r$ defines an equivalence of tensor categories
$\mathcal{M}{}_{\text{num}}(k)\rightarrow\mathcal{M}{}$.
\end{proposition}

\begin{proof}
Such a functor $r$ defines an adequate equivalence relation $\sim$ (see above)
such that $r$ factors into%
\[
\mathcal{M}{}_{\text{rat}}(k)\rightarrow\mathcal{M}{}_{\sim}(k)\overset
{\bar{r}}{\longrightarrow}\mathcal{M}{}%
\]
with $\bar{r}$ a fully faithful tensor functor. Because $P$ is a
pro-reductive, $\mathcal{M}{}$ is semisimple (cf. \cite{deligneM1982}, 2.23).
It follows that $\mathcal{M}{}_{\sim}(k)$ is semisimple (apply the criterion
in \cite{jannsen1992}, Lemma 2), and so $\sim$ is numerical equivalence (ibid.
Theorem 1). The simple objects of $\mathcal{M}{}$ are classified by the orbits
of $\Gal(\mathbb{Q}{}^{\mathrm{al}}/\mathbb{Q}{})$ acting on $X^{\ast}(P)$,
i.e., by the conjugacy classes of Weil $q$-numbers, and so Lemma \ref{0} shows
that $\mathcal{M}{}$ is generated as a tensor category by the images of
abelian varieties. Therefore, $\bar{r}$ is a tensor equivalence, and so
defines an isomorphism of $P$ with $\pi(\mathcal{M}{}_{\text{num}}(k))$. We
can now apply Proposition \ref{00}.
\end{proof}

\begin{remark}
\label{2}When we drop the requirement that $r$ is full, then it is possible to
work with hypotheses much weaker than the full Tate conjecture. Let
$\mathcal{S}{}$ consist of the smooth projective varieties over $\mathbb{F}%
{}_{q}$ whose Frobenius elements are semisimple. In \cite{milne2007rtc} a
notion of a \textquotedblleft good theory of rational Tate classes on
$\mathcal{S}{}$\textquotedblright\ is defined, and it is proved that there
exists at most one such theory. Much of this paper could be rewritten with
\textquotedblleft full Tate conjecture\textquotedblright\ replaced by
\textquotedblleft there exists a good theory rational Tate classes for which
the algebraic classes are rational Tate\textquotedblright\ provided one
removes the requirement that certain functors are full.
\end{remark}

\section{The category of motives}

The next observation goes back to Grothendieck.

\begin{proposition}
\label{02}Let $\mathcal{MM}(\mathbb{F}{}_{q})$ be a pseudo-abelian category
containing $\mathcal{M}_{\text{num}}{}(\mathbb{F}{}_{q})$ as a full
subcategory. Assume

\begin{enumerate}
\item each object $M$ of $\mathcal{\mathcal{M}{}M}(\mathbb{F}{}_{q})$ has a
(weight) filtration%
\[
\cdots\subset W_{i-1}M\subset W_{i}M\subset\cdots
\]
such that $W_{i}M/W_{i-1}M$ is a pure motive of weight $i$;

\item the Frobenius element extends to $\mathcal{M\mathcal{M}{}}{}%
(\mathbb{F}{}_{q})$ and preserves the weight filtrations.
\end{enumerate}

\noindent\noindent Then the inclusion $\mathcal{M}{}_{\text{num}}(\mathbb{F}%
{}_{q})\rightarrow\mathcal{M}{}\mathcal{M}{}(\mathbb{F}{}_{q})$ is an
equivalence of tensor categories.
\end{proposition}

\begin{proof}
For $X$ in $\mathcal{M}{}\mathcal{M}{}(\mathbb{F}{}_{q})$, let $P_{i}(T)$ be
the characteristic polynomial of $\varpi_{W_{i}M/W_{i-1}M}$, and define
$P^{i}(T)$ to satisfy (\ref{e1}). Let $p^{i}=P^{i}(\varpi_{M})$. Then the
$p^{i}$ form a complete set of orthogonal idempotents in $\End(M)$ which
decompose $M$ into a direct sum isomorphic to $\bigoplus\nolimits_{i}%
W_{i}M/W_{i-1}M$.
\end{proof}

\section{Triangulated motivic categories}

Recall that a tensor triangulated category is a category with both a tensor
structure and a triangulated structure satisfying certain compatibilities
(\cite{mazzaVW2006}, Appendix 8A). It is rigid if it admits an internal Hom
or, equivalently, a good theory of duals (\cite{voevodsky2000}, p196). By a
\emph{triangulated motivic category }over a field $k$, we mean a rigid tensor
triangulated category $\mathcal{D}{}$ together with a covariant functor%
\[
R\colon\mathcal{M}{}_{\text{rat}}(k)\rightarrow\mathcal{D}{}%
\]
and isomorphisms for all smooth projective varieties $X$ and all
$i,j\in\mathbb{Z}$%
\begin{equation}
K_{2j-i}(X)^{(j)}\longrightarrow\Hom_{\mathcal{D}}(\1,R(hX)(j)[i])\label{e2}%
\end{equation}
that are natural for the maps defined by algebraic correspondences and reduce
to the identity map when $X$ is a point and $i=j=0$ (see \cite{jannsen2000},
\S 7, p257, which omits the final condition). Here $K_{i}(X)^{(j)}$ is the
subspace of $K_{i}(X)\otimes\mathbb{Q}{}$ on which each Adams operator
$\psi^{m}$ acts as $m^{j}$. According to ibid., p257, over any field $k$ that
admits resolution of singularities, triangulated motivic categories have been
constructed (independently) by \citet{hanamura1995,hanamura1999,hanamura2004},
\citet{levine1998}, and \citet{voevodsky2000}. When $k=\mathbb{F}{}_{q}$,
$\psi^{q}$ acts as $\varpi_{X}$ (\cite{hiller1981}, \S 5; \cite{soule1985},
8.1), and so\footnote{Because the $m^{i}$-eigenspace of $\psi^{m}$ is
independent of $m$ (\cite{seiler1988}, Theorem 1).} $K_{i}(X)^{(j)}$ is the
subspace on which $\varpi_{X}$ acts as $q^{j}$.

Let $\mathcal{\mathcal{D}{}}{}=\mathcal{D}{}(k)$ be a triangulated motivic
category. As we noted in the introduction, for the \textquotedblleft
true\textquotedblright\ triangulated motivic category, there should be a
$t$-structure on $\mathcal{D}{}(k)$ whose heart $\mathcal{M}{}\mathcal{M}%
(k)\overset{\text{{\tiny def}}}{=}\mathcal{D}{}(k)^{\heartsuit}$ is the
category of mixed motives. As Jannsen (2000, \S 7, p257) explains, there
should be the following compatibilities between $R$ and the $t$-structure:

\begin{enumerate}
\item for each standard Weil cohomology, the composite${}$%
\[%
\begin{array}
[c]{lllll}%
\mathcal{M}_{\text{rat}}(k) & \overset{R}{\longrightarrow} & \mathcal{D}{} &
\overset{\bigoplus\nolimits_{i}H^{i}}{\longrightarrow} & \mathcal{M}%
{}\mathcal{M}(k)\\
&  & K & \mapsto & \bigoplus\nolimits_{i}H^{i}(K)
\end{array}
\]
factors through $\mathcal{M}{}_{\text{\textrm{hom}}}(k{})$, and defines a
fully faithful functor $\bar{R}\colon\mathcal{M}{}_{\text{\textrm{hom}}}%
(k{})\rightarrow\mathcal{M}{}\mathcal{M}(k){}$ (here $H^{i}(K)=\tau_{\leq
0}\tau_{\geq0}(K[-i])$);

\item for each smooth projective variety $X$, $\bigoplus\nolimits_{i}{}{}%
H^{i}(R(hX))$ is the weight gradation of $hX$.
\end{enumerate}

\noindent When $k$ is finite, condition (b) says that $H^{i}(R(hX))=\bar
{R}(h^{i}(X))$.

Evidently, there should also be the following compatibilities between the
tensor structures and the $t$-structure:

\begin{enumerate}
\item[(c)] the subcategories $\mathcal{D}{}^{\leq0}$ and $\mathcal{D}{}%
^{\geq0}$ are tensor subcategories of $\mathcal{D}{}$, and $M\mapsto M^{\vee}$
interchanges $\mathcal{D}{}^{\leq0}$ and $\mathcal{D}{}^{\geq0}$, and

\item[(d)] $R\colon\mathcal{M}{}_{\text{rat }}(k)\rightarrow\mathcal{M}%
{}\mathcal{M}{}(k)$ is a tensor functor.
\end{enumerate}

\noindent Note that (c) implies that $\mathcal{\mathcal{M}{}\mathcal{M}{}}(k)$
is a rigid tensor subcategory of $\mathcal{D}{}$.

\begin{definition}
\noindent\label{02m}A $t$-structure on a triangulated motivic category is said
to be \emph{admissible }if it satisfies the conditions (a,b,c,d).
\end{definition}

\begin{theorem}
\label{03}Let $k$ be a finite field. If there exists a triangulated motivic
category $\mathcal{D}{}$ over $k$ and an admissible $t$-structure on
$\mathcal{D}{}$ such that

\begin{itemize}
\item the heart of $\mathcal{D}{}$ is a tannakian category $\mathcal{M}{}$
with fundamental group $P,$ and

\item the functor $\mathcal{M}{}_{\text{rat}}(k)\rightarrow\mathcal{M}%
{}\mathcal{\ }$in (a) above preserves Frobenius elements,
\end{itemize}

\noindent then

\begin{enumerate}
\item the full Tate conjecture holds for all smooth projective varieties over
$k$;

\item for each $l$, the functor $R_{l}\colon\mathcal{M}{}_{\text{hom}%
(l)}(k)\rightarrow\mathcal{M}{}$ defined by $R$ is an equivalence of abelian categories;

\item rational equivalence equals numerical equivalence ($\mathbb{Q}{}$-coefficients);

\item for all $M,N$ in $\mathcal{M(}k)$ and $i\neq0$, $\Hom_{\mathcal{D}{}%
}(M,N[i])=0$.
\end{enumerate}
\end{theorem}

\begin{proof}
Proposition \ref{01} shows that the full Tate conjecture holds and that
$R_{l}$ is essentially surjective (hence an equivalence). Moreover, it allows
us to identify $\mathcal{M}{}$ with $\mathcal{M}_{\text{num}}(k)$.

We next prove (c). When $i=2j$, the isomorphism (\ref{e2}) becomes%
\begin{equation}
K_{0}(X)^{(j)}\simeq\Hom_{\mathcal{D}{}}(\1,R(hX)(j)[2j]). \label{e6}%
\end{equation}
As we noted above, $\varpi_{X}$ acts on $K_{0}(X)^{(j)}$ as $q^{j}$. The Tate
conjecture implies the Lefschetz standard conjecture, and so, for any smooth
projective variety $X$, there exists an isomorphism%
\begin{equation}
R(hX)(j)[2j]\approx\bigoplus\nolimits_{s}h^{s}(X)(j)[2j-s] \label{e6a}%
\end{equation}
(\cite{deligne1968}, \cite{vandenbergh2004}). The characteristic polynomial
$P_{X,s}$ of $\varpi_{X}$ on $h^{s}X$ has roots of absolute value $q^{s/2}$,
and $P_{X,s}(\varpi_{X})$ acts as zero on $h^{s}(X)$ and hence on
$\Hom_{\mathcal{D}{}}(\1,h^{s}(X)(j)[2j-s])$. But we know from (\ref{e6}) that
it acts as $P_{X,s}(q^{j})$. Therefore, $\Hom_{\mathcal{D}{}}(\1,h^{s}%
(X)(j)[2j-s])=0$ unless $s=2j$, and so (\ref{e6}) becomes%
\[
K_{0}(X)^{(j)}\simeq\Hom_{\mathcal{M}_{\text{num}}(k)}(\1,h^{2j}(X)(j)).
\]
Under Grothendieck's isomorphism $K_{0}(X)_{\mathbb{Q}{}}\simeq CH^{\ast
}(X)_{\mathbb{Q}{}}$, the factors $K_{0}(X)^{(j)}$ and $CH^{j}(X)_{\mathbb{Q}%
{}}$ correspond (this is obvious over a finite field), and (by definition)
\[
\Hom_{\mathcal{M}{}_{\text{num}}(k)}(\1,h^{2j}(X)(j))=Z_{\text{num}}%
^{j}(X)_{\mathbb{Q}{}}.
\]
Moreover, our conditions imply that the isomorphism%
\begin{equation}
CH^{j}(X)_{\mathbb{Q}{}}\simeq Z_{\text{num}}^{j}(X)_{\mathbb{Q}{}} \label{e9}%
\end{equation}
obtained by combining these isomorphisms is the canonical one.\footnote{Let
$p$ and $q$ be the projection maps%
\[
\mathrm{pt}\overset{p}{\longleftarrow}X\times\mathrm{pt}\overset
{q}{\longrightarrow}X.
\]
Let $\gamma\in CH^{j}(X)$, and let $f$ be the map $CH^{\ast}(\mathrm{pt}%
)\rightarrow CH^{\ast+j}(X)$ defined by the correspondence $q^{\ast}(\gamma)$.
Then%
\[
f(1_{\mathrm{pt}})\overset{\mathrm{{\scriptsize def}}}{=}q_{\ast}(q^{\ast
}(\gamma)\cup p^{\ast}(\mathrm{1}_{\mathrm{pt}}))=\gamma\cup q_{\ast}p^{\ast
}(1_{\text{pt}})=\gamma\cup1_{X}=\gamma.
\]
As (\ref{e9}) is functorial for correspondences, and the bottom row in%
\[
\begin{CD}
CH^j(X)@>>>Z^j_{\mathrm{num}}(X)\\
@AAfA@AAfA\\
CH^0(\mathrm{pt})@>>>Z^0_{\mathrm{num}}(\mathrm{pt})
\end{CD}
\]
sends $1$ to $1$ (by assumption), it follows that the top row sends $\gamma$
to $\gamma$.
\par
{}} Hence, we have proved (c), and we have shown that%
\begin{equation}
\Hom_{\mathcal{D}{}}(\1,R(hX)(j)[i])=0 \label{e7}%
\end{equation}
when $\ i$ $=2j\neq0$.

Finally, we prove (d). Because $\mathcal{M}{}$ is a rigid subcategory of
$\mathcal{D}{}$, for $M,N$ in $\mathcal{M}{}$ there exists an object
$\underline{\Hom}(M,N)$ in $\mathcal{M}{}$ such that $\Hom_{\mathcal{D}{}%
}(T\otimes M,N)\simeq\Hom_{\mathcal{D}{}}(T,\underline{\Hom}(M,N))$ for all
$T$ in $\mathcal{D}{}$. In particular,%
\[
\Hom_{\mathcal{D}{}}(M,N[i])\simeq\Hom_{\mathcal{D}{}}(\1,\underline
{\Hom}(M,N)[i]).
\]
Therefore, because every object of $\mathcal{M}{}$ is a direct summand of
$R(hX)(j)$ for some smooth projective variety $X$ and integer $j$, it suffices
to prove (d) with $M=\1$ and $N=R(hX)(j)$. We know it when $i=2j$ (see
(\ref{e7})), and so, to complete the proof of (d), it remains to prove
(\ref{e7}) when $i\neq2j$. Because of (\ref{e2}), it suffices to show that (a)
and (c) imply that $K_{i}(X)_{\mathbb{Q}{}}=0$ whenever $i\neq0$. This is done
in \cite{geisser1998}, 3.3. We recall the proof. The functors $K_{i}%
(X)\otimes\mathbb{Q}{}$ factor through $\mathcal{M}{}_{\text{rat}}(k)$
(\cite{soule1984}), and hence (because of (c)) through $\mathcal{M}%
_{\text{num}}(k)$. Therefore, it suffices to prove that $K_{i}(M)\otimes
\mathbb{Q}{}=0$ ($i\neq0$) for $M$ a simple motive in $\mathcal{M}%
_{\text{num}}(k)$. If $M=\mathbb{L}{}^{j}$, then $K_{i}(\mathbb{L}{}^{j})$ is
a direct factor of $K_{i}(\mathbb{P}{}^{j})$, which is torsion
(\cite{quillen1973}). If $M\neq\mathbb{L}{}^{j}$, then $P_{M}(T)$ does not
have $q^{j}$ as a root (\cite{milne1994}, 2.6). As $P_{M}(\varpi_{X})$ acts as
the nonzero rational number $P_{M}(q^{j})$ on $K_{i}(M)^{(j)}$, and also as
zero, the group $K_{i}(M)^{(j)}$ must be zero.
\end{proof}

\begin{corollary}
\label{5}Let $\mathcal{D}$ be as in the theorem, and let $\mathcal{M}{}$ be
its heart. If the inclusion $\mathcal{M}{}\rightarrow\mathcal{D}{}$ extends to
a functor $\mathcal{D}{}^{b}(\mathcal{M}{})\rightarrow\mathcal{D}{}$ (e.g., if
$\mathcal{D}{}$ is endowed with a filtered triangulated category; see \ref{8}a
below), then that functor is an equivalence.
\end{corollary}

\begin{proof}
It suffices to show that $\Hom_{\mathcal{D}{}^{b}(\mathcal{\mathcal{M}{}}%
)}(M,N[i])\rightarrow\Hom_{\mathcal{D}{}}(M,N[i])$ is an isomorphism for all
$M$, $N$ in $\mathcal{M}$ and all $i$ (see \ref{8}b below). For $i=0$ this is
automatic, and for $i\neq0$, both groups are zero (recall that
$\Hom_{\mathcal{D}{}^{b}(\mathcal{M}{})}(M,N[i])\simeq\Ext_{\mathcal{M}{}}%
^{i}(M,N)$, and that $\mathcal{M}{}$ is semisimple).
\end{proof}

\begin{remark}
\label{6}The existence of an admissible $t$-structure on a triangulated
motivic category $\mathcal{D}{}$ implies the existence of a Bloch-Beilinson
filtration on the Chow groups of smooth projective varieties for which%
\begin{equation}
Gr^{s}(CH^{j}(X))\simeq\Hom_{\mathcal{D}{}}(\1,h^{2j-s}(X)(j)[s]) \label{e8}%
\end{equation}
(\cite{jannsen2000}, p258, 4.3). For a finite field, the existence of a
Bloch-Beilinson filtration implies that rational equivalence equals numerical
equivalence ($\mathbb{Q}{}$-coefficients) (ibid., 4.17).
\end{remark}

\begin{remark}
\label{7}Beilinson has conjectured that, for a smooth projective variety $X$,%
\[
Gr^{s}(CH^{j}(X))=\Ext_{\mathcal{M}{}\mathcal{M}_{\text{num}}(k)}%
^{s}(\1,h^{2j-s}(X)(j)).
\]
This is compatible with (\ref{e8}) only if $\mathcal{D}{}=\mathcal{D}{}%
^{b}(\mathcal{M}{}\mathcal{M}_{\text{num}}(k))$ (see the next remark).
\end{remark}

\begin{remark}
\label{8}(a) Let $\mathcal{D}{}$ be a $t$-category with heart $\mathcal{C}$.
Then $D^{b}(\mathcal{C)}$ is also a $t$-category with heart $\mathcal{C}{}$,
but in general there is no obvious relation between $D^{b}(\mathcal{C)}$ and
$\mathcal{D}{}$ (cf. \cite{gelfandM1996}, IV 4.13, p285). In particular, there
will be no obvious functor $r\colon D^{b}(\mathcal{C)}\rightarrow\mathcal{D}%
{}$ extending the inclusion of $\mathcal{C}{}$ into $\mathcal{D}{}$ unless
$\mathcal{D}{}$ is endowed with an additional structure. Beilinson
(1987)\nocite{beilinson1987o} defines the notion of a filtered triangulated
category, and states\footnote{Without proof; cf. the discussion \cite{bbd1982}%
, 3.1, which, however, states that (at that time) the situation had not been
axiomatised.} that such a category over a $t$-category $\mathcal{D}{}$ gives
rise to a well-defined $t$-exact functor $r\colon D^{b}(\mathcal{C)}%
\rightarrow\mathcal{D}$ inducing the identity functor on $\mathcal{C}{}$
(ibid. A.6). The usual triangulated categories are endowed with filtered
triangulated categories over them (ibid. A.2; \cite{bbd1982}, 3.1).

(b) Let $\mathcal{D}$ be a $t$-category with heart $\mathcal{C}{}$. A
$t$-exact functor $r\colon D^{b}(\mathcal{C})\rightarrow\mathcal{D}{}$
inducing the identity functor on $\mathcal{C}$ need not be an equivalence even
when $\mathcal{C}{}$ is semisimple (\cite{deligne1994}, 3.1). We need the
following well-known criterion:%
%TCIMACRO{\TeXButton{bquote}{\begin{quote}}}%
%BeginExpansion
\begin{quote}%
%EndExpansion
Let $r\colon D^{b}(\mathcal{C})\rightarrow\mathcal{D}$ be a $t$-exact functor
inducing the identity functor on $\mathcal{C}$; then $r$ is an equivalence of
$t$-categories if and only if the maps $\Hom_{D^{b}(\mathcal{C})}%
(M,N[i])\rightarrow\Hom_{\mathcal{D}_{\mathcal{C}}^{b}}(M,N[i])$ it defines
are isomorphisms for all $M,N$ in $\mathcal{C}$ and all $i$.%
%TCIMACRO{\TeXButton{equote}{\end{quote}} }%
%BeginExpansion
\end{quote}
%EndExpansion
For $M,N$ in $\mathcal{C}$, let $\Ext_{\mathcal{C}}^{i}(M,N)$ denote the
Yoneda Ext-group, and for $M,N$ in the heart of $\mathcal{D}{}$, let
\[
\Ext_{\mathcal{\mathcal{D}{}}}^{i}(M,N)=\Hom_{\mathcal{D}{}}(M,N[i]).
\]
\noindent Since $\Ext_{\mathcal{C}}^{i}(M,N)\simeq\Hom_{D^{b}(\mathcal{C}%
)}(M,N[i])$ (\cite{verdier1996}, III.3.2.12), the criterion states that
$r\colon D^{b}(\mathcal{C})\rightarrow\mathcal{D}$ is an equivalence of
$t$-categories if and only if the maps $\Ext_{\mathcal{C}}^{i}(M,N)\rightarrow
\Ext_{\mathcal{\mathcal{D}{}}}^{i}(M,N)$\noindent\noindent\ it defines are
isomorphisms for all $M$, $N$, and $i$.
\end{remark}

\section{The motivic $t$-category}

Throughout this section, $k=\mathbb{F}{}_{q}$.

If we want the category of motives to have the Weil-number group $P$ as its
fundamental group, then Corollary \ref{5} shows that $\mathcal{D}{}%
^{b}(\mathcal{M}_{\text{num}}(k))$ is essentially the only candidate for a
triangulated motivic category, and that it will have an admissible
$t$-structure only if the Tate conjecture holds over $k$ and rational
equivalence equals numerical equivalence ($\mathbb{Q}{}$-coefficients). In
this section, we prove that, when we assume these two conjectures,
$\mathcal{D}{}^{b}(\mathcal{M}_{\text{num}}(k))$ does have the hoped for properties.

\begin{proposition}
\label{04}Let $\mathcal{D}{}=D^{b}(\mathcal{M}_{\text{num}}(k))$. Then
$\mathcal{D}{}$ is a  rigid tensor triangulated category with $t$-structure,
and there exists a tensor functor
\[
R\colon{}\mathcal{M}{}_{\text{rat}}(k)\rightarrow\mathcal{D}{},
\]
unique up to a unique isomorphism, such that $H^{i}(RX)=h^{i}(X)[-i]$ for all
$i$.
\end{proposition}

\begin{proof}
Let $C^{b}(\mathcal{M}_{\text{num}}(k))$ be the category of bounded complexes
of objects in $\mathcal{M}_{\text{num}}(k)$, and let $C_{0}^{b}(\mathcal{M}%
_{\text{num}}(k))$ be the full subcategory of bounded complexes whose
differentials are zero. Because $\mathcal{M}_{\text{num}}(k)$ is semisimple,
the functor $D^{b}(\mathcal{M}{}(k))\rightarrow C_{0}^{b}(\mathcal{M}%
_{\text{num}}(k))$ sending $A$ to
\[
\bigoplus\nolimits_{r}H^{r}(A)[-r]=\cdots\rightarrow H^{r-1}(A)\overset
{0}{\longrightarrow}H^{r}(A)\rightarrow\cdots
\]
is an equivalence of categories which is quasi-inverse to the inclusion
functor (\cite{gelfandM1996}, III 2.4, p146). Since $C_{0}^{b}(\mathcal{M}%
_{\text{num}}(k))$ is a direct sum of copies of $\mathcal{M}_{\text{num}}(k)$,
and $\mathcal{M}_{\text{num}}(k)$ is tannakian, it follows that $\mathcal{D}%
{}$ is a rigid tensor category. Define $R$ to be
\begin{equation}
X\mapsto(\cdots\rightarrow h^{r-1}(X)\overset{0}{\longrightarrow}%
h^{r}(X)\rightarrow\cdots). \label{e10}%
\end{equation}
The uniqueness is obvious.
\end{proof}

\begin{remark}
\label{04a}Deligne (1968, 1.11, 1.13)\footnote{This also applies to
$t$-categories. To check this, one only has to check that the spectral
sequence in Deligne's proof exists for $t$-categories (for which there exist
references).} proves the following: \bquote Let $\mathcal{A}{}$ be an abelian
category, and suppose that an object $C$ of $D^{b}(\mathcal{A}{})$ admits
endomorphisms $p_{i}\colon C\rightarrow C$ such that $H^{j}(p_{i})=\delta
_{ij}$ and the $p_{i}$ are orthogonal idempotents; then there is a unique
isomorphism $C\simeq\bigoplus_{i}H^{i}(C)[-i]$ inducing the identity map on
cohomology and such that $p_{i}$ is the $i^{\mathrm{th}}$ projection
map.\equote Let $R^{\prime}$ be a tensor functor $\mathcal{M}{}_{\text{rat}%
}(k)\rightarrow\mathcal{D}{}^{b}(\mathcal{M}{}_{\text{num}}(k))$ and let $R$
be as in (\ref{e10}). Then Deligne's result shows that, for any smooth
projective variety $X$ over $k$, there is a unique isomorphism $R^{\prime
}(X)\simeq R(X)$ inducing the identity on cohomology and such that
$P^{i}(\varpi_{X})$ is the projection from $R^{\prime}(X)$ onto $h^{i}%
(X)[-i]$. Here $P^{i}$ is as in (1).
\end{remark}

\begin{theorem}
\label{06}Assume that the Tate conjecture holds over $k$ and that numerical
equivalence coincides with rational equivalence (with $\mathbb{Q}{}$-coefficients).

\begin{enumerate}
\item $D{}^{b}(\mathcal{M}_{\text{num}}(k))$ has a natural structure of a
triangulated motivic category.

\item The standard $t$-structure on $D{}^{b}(\mathcal{M}_{\text{num}}(k))$ is
admissible (in the sense of \S 4), and it is the unique $t$-structure on
$D{}^{b}(\mathcal{M}_{\text{num}}(k))$ with heart $\mathcal{M}_{\text{num}%
}(k)$.

\item The functor $X\mapsto RX$ sending a smooth projective variety over $k$
to its motivic complex (see \ref{04}) has a unique extension to all varieties
over $k$.

\item For each $l$ (including $p$) there is a $t$-exact functor $R_{l}$ from
$D^{b}(\mathcal{M}_{\text{num}}(k))$ to a $t$-category $\mathcal{D}{}_{l}$
such that $X\mapsto R_{l}(RX)$ is the functor giving rise to the absolute
$l$-adic cohomology.
\end{enumerate}
\end{theorem}

\noindent In the remainder of this section, we explain these statements in
more detail and prove them.

\paragraph{Statement (a).}

We have to construct isomorphisms (\ref{e2}). In computing the right hand side
of (\ref{e2}), we can replace $\mathcal{D}{}^{b}(\mathcal{M}_{\text{num}}(k))$
with the equivalent category $C_{0}^{b}(\mathcal{M}_{\text{num}}%
(k))\simeq\bigoplus\nolimits_{r}\mathcal{M}_{\text{num}}(k)[r]$. Therefore,%
\[
\Hom_{\mathcal{D}{}}(\1,R(X)(j)[i])=\bigoplus\nolimits_{s}\Hom_{\mathcal{D}%
}(\1,h^{s}(X)(j)[i-s]),
\]
and%
\[
\Hom_{\mathcal{D}}(\1,h^{s}(X)(j)[i-s])=\Ext_{\mathcal{M}{}_{\text{num}}%
(k)}^{i-s}(\1,h^{s}(X)(j)).
\]
This group is zero for $i\neq s$ because $\mathcal{M}{}_{\text{num}}(k){}$ is
semisimple, and it is zero for $i=s$, $s\neq2j$, because $\1$ and
$h^{s}(X)(j)$ will then have different weights. It is immediate from the
definition of $\mathcal{M}_{\text{num}}(k)$, that%
\[
\Hom(\1,h^{2j}(X)(j))\simeq Z_{\text{num}}^{j}(X)_{\mathbb{Q}{}}.
\]
On the other hand, $K_{i}(X)_{\mathbb{Q}{}}=0$ for $i\neq0$ (see the proof
\ref{03}), and $K_{0}(X)^{(j)}\simeq CH^{j}(X)_{\mathbb{Q}{}}$. Therefore, we
can define (\ref{e2}) to be the natural map%
\[
CH^{j}(X)_{\mathbb{Q}{}}\rightarrow Z_{\text{num}}^{j}(X)_{\mathbb{Q}{}}%
\]
when $i=2j$ and zero otherwise.

\paragraph{Statement (b).}

By hypothesis, rational, $l$-homological, and numerical equivalence coincide
($\mathbb{Q}{}$-coefficients), and so the standard $t$-structure is obviously
admissible. It is the unique $t$-structure with heart $\mathcal{M}%
{}_{\text{num}}(k)$ because the heart determines the $t$-structure
(\cite{bbd1982}, 1.2, 1.3).

\paragraph{Statement (c).}

We only sketch the argument, leaving the details as an exercise for the
reader. The key point is that \cite{dejong1996}, Theorem 3.1, allows one to
define a simplicial resolution
\[
V\overset{f}{\longleftarrow}U_{\bullet}\overset{j}{\longrightarrow}X_{\bullet}%
\]
of any variety $V$ over $k$ in which $j$ is a simplicial strict
compactification and $f$ is a proper hypercovering of $V$ by a split
simplicial smooth variety (cf. \cite{berthelot1996}, 6.3). One first extends
$R$ to the category of strict compactifications, and then to the simplicial
objects in the category of strict compactifications. Then one defines
$RV=R(U_{\bullet}\rightarrow X_{\bullet})$ for any simplicial resolution
$V\overset{f}{\longleftarrow}U_{\bullet}\overset{j}{\longrightarrow}%
X_{\bullet}$ of $V$. One verifies that $RV\ $is independent of the choice of
the simplicial resolution (up to a well-defined isomorphism), and the map
$V\mapsto RV$ is contravariant for morphisms of varieties.

\paragraph{Statement (d), $l\neq p$.}

For $l\neq p$, let $D_{c}^{b}(k,\mathbb{Z}{}_{l})$ be the category defined in
\cite{deligne1980}, 1.1.2. It is a $t$-category whose heart is the category
$\mathcal{R}{}(k,\mathbb{Z}_{l})$ of finitely generated $\mathbb{Z}{}_{l}%
$-modules endowed with a continuous action of $\Gal(\mathbb{F}{}/k)$. Each
variety $V$ over $k$ defines an object $R\Gamma V$ in $D_{c}^{b}%
(k,\mathbb{Z}{}_{l})$ such that $H^{i}(R\Gamma V)\simeq H_{\mathrm{et}}%
^{i}(V,\mathbb{Z}{}_{l})$ (as objects of $\mathcal{R}{}(k,\mathbb{Z}{}_{l})$).
It is known that $D_{c}^{b}(k,\mathbb{Z}{}_{l})\simeq D^{b}(\mathcal{R}%
{}(k,\mathbb{Z}{}_{l}))$. Now quotient out by the torsion objects to obtain an
equivalence $D_{c}^{b}(k,\mathbb{Q}{}_{l})\simeq D^{b}(\mathcal{R}%
{}(k,\mathbb{Q}_{l}))$ of $\mathbb{Q}{}_{l}$-linear categories. We define
\[
R_{l}\colon D^{b}(\mathcal{M}{}_{\text{\textrm{num}}}(k))\rightarrow
D^{b}(\mathcal{R}{}(k,\mathbb{Q}_{l}))\simeq D_{c}^{b}(k,\mathbb{Q}{}_{l})
\]
to be the derived functor of the fibre functor $\mathcal{M}_{\text{num}}%
{}(k)\rightarrow\mathcal{R}{}(k,\mathbb{Q}{}_{l})$. Applying
\cite{deligne1968}, 1.11, 1.13 (cf. \ref{04a}), we see that, for each smooth
projective variety $V$ over $k$, there is a unique isomorphism $R_{l}%
(V)\simeq\bigoplus_{i}H_{l}^{i}(V)[-i]$ inducing the identity map on
cohomology and such that $P^{i}(\varpi_{V})$ is the $i^{\mathrm{th}}$
projection map. Here $P^{i}$ is the polynomial in (\ref{e1}). This shows that%
\begin{equation}
R\Gamma(V)\simeq R_{l}(RV) \label{e21}%
\end{equation}
when $V$ is projective and smooth. For an arbitrary $V$, we choose a
simplicial resolution $V\overset{f}{\longleftarrow}U_{\bullet}\overset
{j}{\longrightarrow}X_{\bullet}$ of $V$. Because (\ref{e21}) holds for smooth
projective varieties,
\[
R\Gamma(U_{\bullet}\overset{j}{\longrightarrow}X_{\bullet})\simeq
R_{l}(R(U_{\bullet}\overset{j}{\longrightarrow}X_{\bullet})).
\]
Moreover,%
\begin{align*}
R(U_{\bullet}\overset{j}{\longrightarrow}X_{\bullet})  &  \simeq
R(V)\quad\text{(definition of }R(V)\text{)}\\
R\Gamma(U_{\bullet}\overset{j}{\longrightarrow}X_{\bullet})  &  \simeq
R\Gamma(V)\quad\text{(\cite{saintdonat1973}, 4.3.2; also \cite{huber1995},
1.1.3),}%
\end{align*}
and so (\ref{e21}) holds for all varieties.

\paragraph{Statement (d), $l=p$.}

Let $R$ be the Raynaud ring, and $D(R)$ the derived category of the category
of graded $R$-modules (\cite{illusie1983}, 2.1). For a smooth projective
variety $X$ over $k$, let $W\Omega_{X}^{\bullet}$ be the de Rham-Witt complex
on $X$, and let $R\Gamma(W\Omega_{X}^{\bullet})$ be its image under the
derived functor of $\Gamma=\Gamma(X,-)$. Then $R\Gamma(W\Omega_{X}^{\bullet})$
lies in the full subcategory $D_{c}^{b}(R)$ of $D(R)$ consisting of bounded
$R$-complexes whose cohomology modules are coherent (\cite{illusieR1983}, II
2.2), and $H^{i}(R\Gamma(W\Omega_{X}^{\bullet})\simeq H_{\text{crys}}%
^{i}(X/W)$. When we endow $D_{c}^{b}(R)$ with Ekedahl's $t$-structure
(\cite{illusie1983}, 2.4.8) and quotient out by torsion objects, we obtain a
$\mathbb{Q}{}_{p}$-linear $t$-category $D_{c}^{b}(R)_{\mathbb{Q}{}}$ whose
heart is the category $\mathcal{R}(k,\mathbb{\mathbb{Q}{}}{}_{p})$ of
$F$-isocrystals over $k$. It is known that $D_{c}^{b}(R)_{\mathbb{Q}{}}\simeq
D_{\mathcal{R}{}}^{b}(B_{\sigma}[F])$ (derived category of bounded complexes
of $B_{\sigma}[F]$-modules whose cohomology groups are $F$-isocrystals over
$k$; recall $B=W\otimes\mathbb{Q}{}$ and that $B_{\sigma}[F]$ is the twisted
polynomial ring). Define
\[
R_{p}\colon D^{b}(\mathcal{M}{}_{\text{num}}(k))\rightarrow D^{b}%
(\mathcal{R}{}(k,\mathbb{Q}{}_{p}))\rightarrow D_{\mathcal{R}{}}^{b}%
(B_{\sigma}[F])\simeq D_{c}^{b}(R)_{\mathbb{Q}{}}%
\]
to be the composite of the derived functor of the fibre functor $\mathcal{M}%
{}_{\text{num}}(k)\rightarrow\mathcal{R}{}(k,\mathbb{Q}{}_{p})$ with the
natural functors. The proof can now be completed as in the case $l\neq p$
except that the reference to \cite{saintdonat1973} must be replaced by a
reference to \cite{tsuzuki2003}.

\begin{remark}
\label{9m}Statement (c) and (d) of the theorem are very strong. Consider, for
example, a closed subvariety $Z$ of codimension $r$ in a smooth projective
variety $X$ and a smooth alteration $\tilde{Z}\rightarrow Z$. Then the theorem
says that there is an exact sequence%
\[
h^{i-2r}(\tilde{Z})(r)\rightarrow h^{i}(X)\rightarrow h^{i}(U),\quad
U=X\smallsetminus Z,
\]
whose $l$-adic realization is the sequence in (\ref{gtc2a}) for $l\neq p$.
\end{remark}

\subsubsection{Application.}

\begin{plain}
\label{9n}Using (c) and (d), we can extend the definition of $\mathbb{Q}{}%
_{p}$ cohomology (\cite{milne1986}, p309) from smooth projective varieties to
all varieties, namely, for any variety $X$ over $k$, define%
\[
H^{i}(X,\mathbb{Q}_{p}(r))=\Hom_{\mathcal{D}{}(k;\mathbb{Q}{}_{p})}%
(\1,R_{p}(RX)(r)[i]).
\]
The main theorem of \cite{milneR2005} shows that this agrees with the original
definition when $X$ is smooth and projective.
\end{plain}

\section{The $\mathbb{Q}{}$-algebra of correspondences at the generic point}

In this section, $k=\mathbb{F}_{q}$ and we assume that the Tate conjecture
holds over $k$ and that numerical equivalence equals rational equivalence
($\mathbb{Q}{}$-coefficients). We allow $l=p$.

\subsection{Effective motives}

Let $\mathcal{M}{}^{\mathrm{eff}}(k)$ be the category of effective motives
given by Grothendieck's construction using algebraic classes modulo numerical
equivalence as correspondences. It is an abelian nonrigid tensor category, and
we let $\mathcal{D}{}^{\text{eff}}(k)=D{}^{b}(\mathcal{M}{}^{\text{eff}}(k))$.
Much of Theorem \ref{06} continues to hold. In particular, attached to a
smooth projective variety $X$ and an open subvariety $U$, there is a
well-defined restriction map $h^{i}(X)\rightarrow h^{i}(U)$ whose $l$-adic
realization is $H_{l}^{i}(X)\rightarrow H_{l}^{i}(U)$ (cf. \ref{9m},
\ref{9n}). We define
\[
F_{a}^{r}h^{i}(X)=\bigcup\nolimits_{U}\Ker(h^{i}(X)\rightarrow h^{i}(U))
\]
where $U$ runs over the open subvarieties of $X$ such that $X\smallsetminus U$
is of codimension at least $r$.

\begin{proposition}
\label{10}For all $l$ (including $l=p$)%
\[
R_{l}(F_{a}^{r}h^{i}(X))=F_{b}^{r}H_{l}^{i}(X).
\]

\end{proposition}

\begin{proof}
The functor $R_{l}$ is exact, and so
\[
R_{l}(F_{a}^{r}h^{i}(X))=F_{a}^{r}H_{l}^{i}(X).
\]
Therefore, the statement follows from the generalized Tate conjecture
(\ref{gtc8}, \ref{gtc2a}, \ref{gtc9c}).
\end{proof}

\begin{corollary}
\label{10a}For any smooth projective variety $X$ over $k$, $F_{a}^{r}h^{i}(X)$
is the largest effective submotive of $h^{i}(X)$ of the form $M(-r)$ for some
effective motive $M$.
\end{corollary}

\begin{proof}
Obvious.
\end{proof}

\subsection{Definition of the $\mathbb{Q}{}$-algebra of correspondences at the
generic point}

In this subsection, we translate some definitions and results of
\cite{beilinson2002} into our context. Let $X$ be a connected algebraic
variety of dimension $n$ over a finite field $k$, and let $\eta$ be its
generic point. Define%
\[
CH^{n}(\eta\times\eta)=\varinjlim CH^{n}(U\times U),
\]
where $U$ runs over the open subvarieties of $X$. Following
\cite{beilinson2002}, 1.4, we define
\[
A(X)=CH^{n}(\eta\times\eta)\otimes\mathbb{Q}{}.
\]
Composition of correspondences makes $A(X)$ into an associative $\mathbb{Q}%
$-algebra, called the $\mathbb{Q}$\emph{-algebra of correspondences at the
generic point. }

Denote by $\bar{h}^{n}(X)$ the image of the canonical map $h^{n}(X)\rightarrow
h^{n}(\eta)$ (ind object of $\mathcal{M}{}^{\mathrm{eff}}(k{})$).

\begin{theorem}
\label{11}For any connected smooth projective varieties $X,X^{\prime}$ of
dimension $n$ over $k$, the map%
\[
CH^{n}(\eta^{\prime}\times\eta)\otimes\mathbb{Q}{}\rightarrow\Hom(\bar{h}%
^{n}(\eta),\bar{h}^{n}(\eta^{\prime}))
\]
is an isomorphism.
\end{theorem}

\begin{proof}
Beilinson's proof (2002, 4.9) applies in our context.
\end{proof}

\begin{corollary}
\label{12}For any connected smooth projective variety $X$ of dimension $n$
over $k$, there is a canonical isomorphism of $\mathbb{Q}{}$-algebras%
\[
A(X)\simeq\End(\bar{h}^{n}(X)).
\]

\end{corollary}

\begin{proof}
It is only necessary to observe that composition of correspondences
corresponds to composition of endomorphisms (\cite{beilinson2002}, 4.10).
\end{proof}

\begin{corollary}
\label{12m}The $\mathbb{Q}{}$-algebra $A(X)$ is finite-dimensional and semisimple.
\end{corollary}

\begin{proof}
Immediate from (\ref{12}) because $\mathcal{M}{}^{\mathrm{eff}}(k)$ is a
semisimple category over $\mathbb{Q}{}$ with finite-dimensional $\Hom$s.
\end{proof}

\subsection{Calculation of the $\mathbb{Q}{}$-algebra of correspondences at
the generic point}

\begin{proposition}
\label{13}For a connected curve $X$ over $k$,
\[
A(X)\simeq\End(J)\otimes\mathbb{Q}{}%
\]
where $J$ is the Jacobian of a smooth complete model of $X$.
\end{proposition}

\begin{proof}
As $X$ is geometrically reduced, its smooth locus $X^{\prime}$ can be embedded
in a smooth projective curve $Y$, and $X^{\prime}\hookrightarrow Y$ is
uniquely determined up to a unique isomorphism. As%
\[
A(X)\simeq A(X^{\prime})\simeq A(Y)
\]
we may as well assume that $X$ itself is smooth and projective. For any
nonempty open $U$, the map $h^{1}(X)\rightarrow h^{1}(U)$ is injective because
$H_{l}^{1}(X)\rightarrow H_{l}^{1}(U)$ is injective, and so $\bar{h}%
^{1}(X)=h^{1}(X)$. Therefore, $A(X)\simeq\End(h^{1}(X))$, and it follows from
the isomorphism%
\[
CH^{1}(X\times X)\simeq CH^{1}(X)\oplus CH^{1}(X)\oplus\End(J)
\]
(\cite{weil1948s}), that%
\[
\End(h^{1}(X))\simeq\End(J)\otimes\mathbb{Q}{}.
\]

\end{proof}

For a connected smooth projective variety $X$ of dimension $n$ over $k$,
define\footnote{For $l\neq p$, this is $Gr^{0}H_{l}^{n}(X)$, the
\textquotedblleft composante pure de niveau $n$\textquotedblright\ of
$H_{l}^{n}(X),$ of \citet[p162]{grothendieck1968}.}%
\[
\bar{H}_{l}^{n}(X)=H_{l}^{n}(X)/F_{b}^{0}H_{l}^{n}(X).
\]
For $l\neq p$, the quotient map $H_{l}^{n}(X)\rightarrow\bar{H}_{l}^{n}(X)$
defines an isomorphism of $\bar{H}_{l}^{n}(X)$ with the Tate substructure of
$H_{l}^{n}(X)$ whose Frobenius eigenvalues $\alpha$ are such that $a/q$ is not
an algebraic integer. The quotient map $H_{p}^{n}(X)\rightarrow\bar{H}_{p}%
^{n}(X)$ can be identified with the map
\[
H^{n}(X,W\Omega^{\bullet})_{\mathbb{Q}{}}\rightarrow H^{n}(X,W\mathcal{O}%
{}_{X})_{\mathbb{Q}{}}\simeq H_{p}^{n}(X)_{[0,1[}%
\]
(\cite{illusie1979}, II 3.5.3, p616).

\begin{proposition}
\label{14}For all primes $l$ (including $l=p$),%
\[
R_{l}(\bar{h}^{n}(X))\simeq\bar{H}_{l}^{n}(X).
\]

\end{proposition}

\begin{proof}
Clearly,%
\[
0\rightarrow F_{a}^{0}h^{n}(X)\rightarrow h^{n}(X)\rightarrow\bar{h}%
^{n}(X)\rightarrow0
\]
is exact. On applying the exact functor $R_{l}$, this gives an exact sequence%
\[
0\rightarrow F_{b}^{0}H_{l}^{n}(X)\rightarrow H_{l}^{n}(X)\rightarrow\bar
{H}_{l}^{n}(X)\rightarrow0
\]
by (\ref{12}).
\end{proof}

\begin{theorem}
\label{15}For all primes $l$ (including $l=p$),%
\[
A(X)\otimes\mathbb{Q}{}_{l}\simeq\End(\bar{H}_{l}^{n}(X))
\]
(endomorphisms of $\bar{H}_{l}^{n}(X)$ as a Tate structure when $l\neq p$;
endomorphisms of $\bar{H}_{p}^{n}(X)$ as an $F$-isocrystal when $l=p$).
\end{theorem}

\begin{proof}
Follows from Proposition \ref{14} and the fact that $R_{l}$ defines
isomorphisms%
\[
\Hom(M,N)\otimes_{\mathbb{Q}{}}\mathbb{Q}{}_{l}\simeq\Hom(R_{l}M,R_{l}N).
\]

\end{proof}

\begin{example}
\label{16}If $H^{n}(X,W\mathcal{O}{}_{X})$ is torsion, then $A(X)=0$. This is
the case, for example, if $X$ is a supersingular abelian surface, a
supersingular $K3$ surface, or an Enriques surface (\cite{illusie1979}, 7.1,
7.2, 7.3).
\end{example}

\begin{remark}
\label{16m}It is possible to recover the rank of a motive $M$ from its
endomorphism algebra $\End(M)$. According to the Wedderburn theorems,
\[
\End(M)=\prod\nolimits_{j}M_{r_{j}}(D_{j})
\]
with each $D_{j}$ a division algebra over $\mathbb{Q}{}$. If $Z_{j}$ is the
centre of $D_{j}$, then%
\[
\rank(M)=\sum\nolimits_{j}r_{j}\cdot\lbrack Z_{j}\colon\mathbb{Q}]\cdot\lbrack
D_{j}\colon Z_{j}]^{1/2}.
\]

\end{remark}

\begin{remark}
\label{16n}Since $A(X)$ is a birational invariant, (\ref{12}) and\ (\ref{16m})
show that the rank of $\bar{h}^{n}(X)$ ($n=\dim X$) is a birational invariant
of connected smooth projective varieties. Hence the same is true of its
$p$-adic realization, i.e.,%
\[
\rank H^{n}(X,W\mathcal{O}{}_{X})=\rank H^{0}(X,W\Omega_{X}^{n})
\]
is a birational invariant of connected smooth projective varieties over a
finite field. Of course, it is classical that%
\[
\dim H^{n}(X,\mathcal{O}{}_{X})=\dim H^{0}(X,\Omega_{X}^{n})
\]
is a birational invariant (\cite{hartshorne1977}, II Ex 8.8), but%
\[
\rank H^{n}(X,W\mathcal{O}{}_{X})\neq\dim H^{n}(X,\mathcal{O}{}_{X}),
\]
for example, when $n=2$ and $X$ is a supersingular abelian surface.
\nocite{illusie1979}Illusie (1979, II 2.18, p614) proves that $H^{0}%
(X,W\Omega_{X}^{n})$ is of finite-type over $W$ with $F$ acting as an
automorphism. The formal $p$-divisible group $G$ with Cartier module
$H^{0}(X,W\Omega_{X}^{n})/\mathrm{torsion}$ has%
\[
\dim(G)=\rank(G)=\rank H^{0}(X,W\Omega_{X}^{n})
\]
(cf. ibid. II 4.4, p621) and so $\dim(G)$ and $\rank(G)$ are also birational invariants.
\end{remark}

\paragraph{Explicit description of $A(X)$}

\begin{plain}
Let $X$ be a smooth projective variety over $k$, and let $(\alpha_{i})_{1\leq
i\leq\beta_{n}}$ be the family of eigenvalues of $\varpi_{X}$ on $H_{l}%
^{n}(X)$. Then the family $S(X)$ of eigenvalues of $\varpi_{\bar{h}^{n}(X)}$
consists of the $\alpha_{i}$ for which $\alpha_{i}/q$ is not an algebraic
integer. Therefore, by \cite{milne1994}, 2.14--2.15, the semisimple
$\mathbb{Q}{}$-algebra $A(X)=\End(\bar{h}(X))$ has the following description.
Let $o_{1},\ldots,o_{s}$ be the distinct orbits for the action of
$\Gal(\mathbb{Q}{}^{\mathrm{al}}/\mathbb{Q}{})$ on $S(X)$ and let $r_{j}$ be
the multiplicity of $o_{j}$:
\[
F(X)=\coprod\nolimits_{j}r_{j}o_{j}.
\]
Then
\[
\bar{h}^{n}(X)=\sum\nolimits_{j}r_{j}N_{j}%
\]
where $N_{j}$ is a simple motive with Frobenius eigenvalues the elements of
$o_{j}$, and%
\[
A(X)\simeq\prod\nolimits_{j}M_{r_{j}}(\End(N_{j})).
\]
Let $\alpha\in o_{j}$. Then $\End(N_{j})$ is isomorphic to a central simple
algebra $D_{j}$ over $Z_{j}=\mathbb{Q}{}[\alpha]$ with invariants (at the
primes $v$ of $\mathbb{Q}{}[\alpha]$)%
\[
\inv_{v}(D_{j})=\left\{
\begin{array}
[c]{lll}%
\frac{1}{2} &  & \text{if }v\text{ is real and }n\text{ is odd}\\
\frac{\ord_{v}(\alpha)}{\ord_{v}(q)}\cdot\lbrack\mathbb{Q}{}[\alpha
]_{v}:\mathbb{Q}{}_{p}] &  & \text{if }v|p\\
0 &  & \text{otherwise.}%
\end{array}
\right.
\]
Therefore, the degree $[\mathbb{Q}{}[\alpha]\colon\mathbb{Q}{}]$ is the order
of $o_{i}$, and the degree $[D_{j}\colon\mathbb{Q}{}[\alpha]]=e^{2}$ where $e$
is the least common denominator of the numbers $\inv_{v}(D_{j})$.
\end{plain}

Following Beilinson (2002, p37), the pessimists will be tempted to look for
counter-examples to the above calculations in order to ruin the conjectures.

\section{Base fields algebraic over a finite field}

Let $k$ be a subfield of $\mathbb{F}{}$, and assume that the Tate conjecture
holds and numerical equivalence equals rational equivalence ($\mathbb{Q}{}%
$-coefficients) for finite subfields of $k$. When we define the various
categories for $k$ to be the $2$-category direct limits of the categories for
$k^{\prime}$ with $k^{\prime}$ running over the finite subfields of $k$, then
these categories for $k$ inherit the properties of the corresponding
categories for $k^{\prime}$.

\bsmall
\bibliographystyle{cbe}
\bibliography{/Current/refs}

\bigskip\noindent James S. Milne, \newline Mathematics Department, University
of Michigan, Ann Arbor, MI 48109, USA, \newline Email: \url{jmilne@umich.edu}
\newline Webpage: \url{www.jmilne.org/math/}

\medskip\noindent Niranjan Ramachandran, \newline Mathematics Department,
University of Maryland, College Park, MD 20742, USA, \newline Email:
\url{atma@math.umd.edu}, \newline Webpage: \url{www.math.umd.edu/~atma/}

\esmall
\end{document}